\documentclass[10pt]{amsart}
\usepackage{amsfonts}
\usepackage{ifthen}
\usepackage{amsthm}
\usepackage{amsmath}
\usepackage{graphicx}
\usepackage{amscd,amssymb,amsthm}
\usepackage{graphicx}
\usepackage{epstopdf}
\usepackage{hyperref}
\usepackage{color}
\newcounter{alphabet}
\newcounter{tmp}

\newcounter{minutes}
\divide\time by 60
\newcounter{hours}
\multiply\time by 60 \addtocounter{minutes}{-\time}

\newtheorem{theorem}{Theorem}
\newtheorem{corollary}{Corollary}

\newcommand{\real}{\operatorname{Re}}

\newcommand{\lp}{\mathcal{LP}}


\title{Geometric and monotonic properties of hyper-Bessel functions}

\author[\.{I}. Akta\c{s}]{\.{I}brah\.{I}m Akta\c{s}}
\address{Department of Mathematical Engineering, Faculty of Engineering and Natural Sciences, G\"{u}m\"{u}\c{s}hane University, G\"{u}m\"{u}\c{s}hane, Turkey}
\email{aktasibrahim38@gmail.com}

\author[\'A. Baricz]{\'Arp\'ad Baricz}
\address{Department of Economics, Babe\c{s}-Bolyai University, Cluj-Napoca, Romania}
\address{Institute of Applied Mathematics, \'Obuda University, Budapest, Hungary}
\email{bariczocsi@yahoo.com}

\author[S. Singh]{Sanjeev Singh}
\address{Discipline of Mathematics, Indian Institute of Technology Indore,
Indore, India}
\email{snjvsngh@iiti.ac.in}

\begin{document}

\maketitle

\begin{abstract}
Some geometric properties of a normalized hyper-Bessel functions are investigated. Especially we focus on the radii of starlikeness, convexity, and uniform convexity of hyper-Bessel functions and we show that the obtained radii satisfy some transcendental equations. In addition, we give some bounds for the first positive zero of normalized hyper-Bessel functions, Redheffer-type inequalities, and bounds for this function. In this study we take advantage of Euler-Rayleigh inequalities and Laguerre-P\'{o}lya class of real entire functions, intensively.

\end{abstract}

\section{\bf Introduction}
Let $\mathbb{D}_r$ be the open disk $\{z\in\mathbb{C}:\left|z\right|<r\}$ with radius $r>0$ and $\mathbb{D}_1=\mathbb{D}$. Let $\mathcal{A}$ denote the class of analytic functions $f:\mathbb{D}_r\rightarrow\mathbb{C},$
which satisfy the normalization conditions $f(0)=f^{\prime}(0)-1=0$. By $\mathcal{S}$ we mean the class of functions belonging to $\mathcal{A}$, which are univalent in $\mathbb{D}_r$ and let $\mathcal{S}^*$ be the subclass of $\mathcal{S}$ consisting of functions which are starlike with respect to the origin in $\mathbb{D}_r$. The analytic characterization of the class of starlike functions is
$$
\mathcal{S}^*=\left\{f\in\mathcal{S}:\real\left(\frac{zf^{\prime}(z)}{f(z)}\right)>0\text{ for all }z\in\mathbb{D}_r\right\}.
$$
The real number
$$
r^{\ast }(f)=\sup \left\{ r>0 :\real\left(\frac{zf^{\prime }(z)}{f(z)}\right) >0 \;\text{for all }z\in\mathbb{D}_r\right.\}
$$
is called the radius of starlikeness of the function $f$. Note that $r^{\ast }(f)$ is the largest radius such that the image region $f(\mathbb{D}_{r^*(f)})$ is a starlike domain with respect to the origin. On the other hand, the class of convex functions is defined by
$$
\mathcal{K}=\bigg\{f\in \mathcal{S}:\real\left(1+\frac{zf^{\prime \prime}(z)}{f^{\prime}(z)}\right)>0 \text{ for all } z\in \mathbb{D}_r \bigg\}.
$$
It is known that the convex functions do not need to be normalized, namely, the definition of $\mathcal{K}$ is also valid for non-normalized analytic function $f:\mathbb{D}_r\rightarrow\mathbb{C}$ which has the property $f^{\prime}(0)\neq0$. The radius of convexity of an analytic locally univalent function $f:\mathbb{C}\rightarrow\mathbb{C}$ is defined by
$$
r^{c}(f)=\sup\bigg\{r>0:\real\left(1+\frac{zf^{\prime \prime}(z)}{f^{\prime}(z)}\right)>0 \text{ for all } z\in \mathbb{D}_r \bigg\}.
$$
Note that $r^{c}(f)$ is the largest radius for which the image domain $f\left(\mathbb{D}_{r^{c}(f)}\right)$ is a convex domain in $\mathbb{C}.$ For more information about starlike and convex functions we refer to Duren's book \cite{Duren} and to the references therein.

Goodman \cite{Goodman} has introduced the concept of uniform convexity for the functions in class $\mathcal{A}$ . A function $f$ is said to be uniformly convex in $\mathbb{D}_r$ if $f$ is in the class of usual convex functions and has the property that for every circular arc $\gamma$ contained in $\mathbb{D}_r$, with the center $\zeta$ also in $\mathbb{D}_r$, the arc $f(\gamma)$ is a convex arc.
An analytic description of the uniformly convex functions has been given by R\o{}nning \cite{Ronning} and it reads as follows: the function $f$  is a uniformly convex functions if and only if
\begin{equation*}\label{1.2}
\real\left(1+\frac{zf^{\prime\prime}(z)}{f^{\prime}(z)}\right)>\left|\frac{zf^{\prime\prime}(z)}{f^{\prime}(z)}\right| \text{ for all } z\in\mathbb{D}_r.
\end{equation*}
The radius of uniform convexity of a function $f$ is defined by (see for example \cite{Deniz})
$$r^{uc}(f)=\sup\Big\{r>0: \real\left(1+\frac{zf^{\prime\prime}(z)}{f^{\prime}(z)}\right)>\left|\frac{zf^{\prime\prime}(z)}{f^{\prime}(z)}\right| \text{ for all } z\in\mathbb{D}_r\Big\}.$$

Recently, there has been a vivid interest on some geometric properties such as univalency, starlikeness, convexity and uniform convexity of various special functions such as Bessel, Struve, Lommel, Wright and $ q $-Bessel functions (see \cite{aktas3,aktas2,aktas1,aktas4,mathematica,publ,lecture,bdoy,bsk,btk,bos,samy,basz}). In the above mentioned papers  the authors have used frequently some properties of the zeros of these functions. In this paper our aim is to study some properties of hyper-Bessel functions. The location of the zeros of hyper-Bessel functions was studied by Chaggara and Romdhane \cite{Chaggara}. Motivated by the earlier works, our purpose is to see how far we can go with the extensions of the results on the Bessel functions to hyper-Bessel functions. More precisely, we obtain some new results on the radii of starlikeness, convexity and uniform convexity for the normalized  hyper-Bessel functions. Also, we present new bounds for the first positive zero of normalized hyper-Bessel function,  interlacing properties of the zeros of normalized hyper-Bessel functions, Redheffer-type inequalities, and bounds for this function. The results presented in this paper are natural extensions of known results for the classical Bessel functions. The conclusion of the paper is that many of the properties of classical Bessel functions can be extended in a natural way to hyper-Bessel functions and the reason why we can do this is that the hyper-Bessel functions have also the beautiful infinite sum and product representation structure. The paper is organized as follows: this section contains some information about the hyper-Bessel functions and their zeros, section 2 is devoted for the main results of the paper, while the third section contains the proofs of the main results.

Now, consider the hyper-Bessel function defined by \cite{Chaggara}
\begin{equation}\label{eq1.1}
J_{\alpha_d}(z)=\dfrac{\left(\frac{z}{d+1}\right)^{\alpha_1+\dots+\alpha_d}}{\Gamma\left(\alpha_1+1\right)\cdots\Gamma\left(\alpha_d+1\right)}{_0F_d}
\left(\begin{matrix}-\\(\alpha_{d}+1)\end{matrix};-\left(\frac{z}{d+1}\right)^{d+1}\right),
\end{equation}
where the notation
\begin{equation}\label{eq1.2}
{_pF_q}\left(\begin{matrix}(\beta_p)\\(\gamma_q)
\end{matrix};x\right)=\sum_{n\geq0}\frac{(\beta_1)_n(\beta_2)_n\cdots(\beta_p)_n}{(\gamma_1)_n(\gamma_2)_n\cdots(\gamma_p)_n}\frac{x^n}{n!}
\end{equation}
denotes the generalized hypergeometric function,
$(\beta)_n$ is the shifted factorial
(or Pochhammer's symbol) defined by $(\beta)_0=1, (\beta)_n=\beta(\beta+1)\cdots(\beta+n-1), n\geq1$ and the contracted notation $\alpha_d$ is used to abbreviate the array of $d$ parameters $\alpha_1,\dots, \alpha_d$.

By using the equations \eqref{eq1.1} and \eqref{eq1.2} one can obtained that the function $z\mapsto J_{\alpha_d}(z)$ has the following infinite sum representation
\begin{equation}\label{eq 1.3}
J_{\alpha_d}(z)=\sum_{n\geq0}\frac{(-1)^n}{n!\Gamma\left(\alpha_1+1+n\right)\cdots\Gamma\left(\alpha_d+1+n\right)}\left(\frac{z}{d+1}\right)^{n(d+1)+\alpha_1+\dots+\alpha_d}.
\end{equation}
By choosing $d=1$ and putting $\alpha_1=\nu$ in the expression \eqref{eq 1.3} we get the classical Bessel function
\begin{equation*}\label{eq 1.4}
J_\nu(z)=\sum_{n\geq0}\frac{(-1)^n}{n!\Gamma\left(\nu+n+1\right)}\left(\frac{z}{2}\right)^{2n+\nu}.
\end{equation*}
The normalized hyper-Bessel function $\mathcal{J}_{\alpha_d}(z)$ is defined by
\begin{equation}\label{eq 1.5}
J_{\alpha_d}(z)=\frac{\left(\frac{z}{d+1}\right)^{\alpha_1+\dots+\alpha_d}}{\Gamma\left(\alpha_1+1\right)\cdots\Gamma\left(\alpha_d+1\right)}\mathcal{J}_{\alpha_d}(z).
\end{equation}
By combining the equations \eqref{eq 1.3} and \eqref{eq 1.5} we obtain the following infinite sum representation
\begin{equation}\label{eq1.6}
\mathcal{J}_{\alpha_d}(z)=\sum_{n\geq0}\frac{(-1)^n}{n!(\alpha_1+1)_n\cdots(\alpha_d+1)_n}\left(\frac{z}{d+1}\right)^{n(d+1)}.
\end{equation}
Since the function $ \mathcal{J}_{\alpha_d} $ does not belong to the class $\mathcal{A}$, we consider the following normalized form
\begin{equation}\label{eq1.7}
{f_{\alpha_d}(z)}=z\mathcal{J}_{\alpha_d}(z)=\sum_{n\geq0}\frac{(-1)^n}{n!\left(d+1\right)^{n(d+1)}(\alpha_1+1)_n\cdots(\alpha_d+1)_n}z^{n(d+1)+1}
\end{equation}
so that the function ${f_{\alpha_d}}\in\mathcal{A}$.

The function $J_{\alpha_d}$ has the infinite product representation \cite[Equation(5.5)]{Chaggara}
\begin{equation}\label{eq 1.8}
J_{\alpha_d}(z)=\frac{\left(\frac{z}{d+1}\right)^{\alpha_1+\dots+\alpha_d}}{\Gamma\left(\alpha_1+1\right)\cdots\Gamma\left(\alpha_d+1\right)}\prod_{n\geq1}\left(1-\frac{z^{d+1}}{j_{\alpha_d,n}^{d+1}}\right),
\end{equation}
where $j_{\alpha_d,n}$ denotes the $n$th positive zero of the function $\mathcal{J}_{\alpha_d}$.

In view of the equations \eqref{eq 1.5} and \eqref{eq 1.8}, we have
\begin{equation}\label{eq 1.9}
\mathcal{J}_{\alpha_d}(z)=\prod_{n\geq1}\left(1-\frac{z^{d+1}}{j_{\alpha_d,n}^{d+1}}\right)
\end{equation}
and consequently
\begin{equation}\label{eq 1.10}
{f_{\alpha_d}}(z)=z\prod_{n\geq1}\left(1-\frac{z^{d+1}}{j_{\alpha_d,n}^{d+1}}\right).
\end{equation}

Now, we would like to mention that by using the equations \eqref{eq1.7} and \eqref{eq 1.10} we can obtain the following Euler-Rayleigh sums for the positive zeros of the function ${f_{\alpha_d}}$. From the equality \eqref{eq1.7} we have
\begin{equation}\label{eq.1.11}
{f_{\alpha_d}}(z)=z-\frac{z^{d+2}}{(d+1)^{d+1}(\alpha_1+1)_1\cdots(\alpha_d+1)_1}+\frac{z^{2d+3}}{2(d+1)^{2(d+1)}(\alpha_1+1)_2\cdots(\alpha_d+1)_2}-\cdots
\end{equation}
Now, if we consider \eqref{eq 1.10}, then some calculations yield that
\begin{equation}\label{eq.1.12}
{f_{\alpha_d}}(z)=z-\sum_{n\geq1}\frac{1}{{j_{\alpha_d,n}^{d+1}}}z^{d+2}+\frac{1}{2}\left(\left(\sum_{n\geq1}\frac{1}{{j_{\alpha_d,n}^{d+1}}}\right)^2-\sum_{n\geq1}\frac{1}{{j_{\alpha_d,n}^{2(d+1)}}}\right)z^{2d+3}-\cdots
\end{equation}
By equating the first few coefficients with the same degrees in equations \eqref{eq.1.11} and \eqref{eq.1.12} we get,
\begin{equation}\label{eq.1.13}
\sum_{n\geq1}\frac{1}{{j_{\alpha_d,n}^{d+1}}}=\frac{1}{(d+1)^{d+1}(\alpha_1+1)\cdots(\alpha_d+1)}
\end{equation}
and
\begin{equation}\label{eq.1.14}
\frac{1}{2}\left(\left(\sum_{n\geq1}\frac{1}{{j_{\alpha_d,n}^{d+1}}}\right)^2-\sum_{n\geq1}\frac{1}{{j_{\alpha_d,n}^{2(d+1)}}}\right)=\frac{1}{2(d+1)^{2(d+1)}(\alpha_1+1)_2\cdots(\alpha_d+1)_2}.
\end{equation}
respectively.
Finally, substituting the \eqref{eq.1.13} in the \eqref{eq.1.14} we have
\begin{equation}\label{eq.1.15}
\sum_{n\geq1}\frac{1}{{j_{\alpha_d,n}^{2(d+1)}}}=\frac{B(\alpha_d)-A(\alpha_d)}{(d+1)^{2(d+1)}\left(A(\alpha_d)\right)^2B(\alpha_d)},
\end{equation}
 where
 $$
 A(\alpha_d)=\prod_{i=1}^{d}\left(\alpha_{i}+1\right)\text{ and }B(\alpha_d)=\prod_{i=1}^{d}\left(\alpha_{i}+2\right).
 $$
Here, it is important mentioning that for $d=1$ and $\alpha_1=\nu$ the equations \eqref{eq.1.13} and \eqref{eq.1.15} reduce to
$$
\sum_{n\geq1}\frac{1}{{j_{\nu,n}^{2}}}=\frac{1}{4(\nu+1)}
$$
and
$$
\sum_{n\geq1}\frac{1}{{j_{\nu,n}^{4}}}=\frac{1}{16(\nu+1)^2(\nu+2)}
$$
respectively, where $j_{\nu,n}$ denotes the $n$th zero of classical Bessel function $J_\nu.$

Now let us recall the Laguerre-P\'{o}lya class $\lp$ of real entire functions which will be used
in the sequel . A real entire function $q$ belongs to the Laguerre-P\'{o}lya class $\lp$ if it can be represented in the form
$$
q(z) = c z^{m} e^{-a z^{2} + \beta z} \prod_{k\geq1}
\left(1+\frac{z}{z_{k}}\right) e^{-\frac{z}{z_{k}}},
$$
with $c,$ $\beta,$ $z_{k} \in \mathbb{R},$ $a \geq 0,$ $m\in\{0,1,2,\dots\}
,$ $\sum\limits_{k\geq1} z_{k}^{-2} < \infty.$

\section{\bf Properties of hyper-Bessel functions}
\setcounter{equation}{0}

\subsection{Radii of starlikeness of normalized hyper-Bessel functions}
Our first main result is about the radius of starlikeness $r^{\star}\left(f_{\alpha_d}\right)$ of the normalized hyper-Bessel function $z\mapsto{f_{\alpha_d}(z)}$.

\begin{theorem}\label{th1}
Let $\alpha_{i}>-1$ for $i\in \{1,2,\dots,d\}.$ Then the radius of starlikeness $r^{\star}\left(f_{\alpha_d}\right)$ of the normalized hyper-Bessel function
$z\mapsto{f_{\alpha_d}(z)}=z\mathcal{J}_{\alpha_d}(z)$ is the smallest positive root of the equation $$z\mathcal{J}_{\alpha_d}^{\prime}(z)+\mathcal{J}_{\alpha_d}(z)=0$$
and satisfies the following inequalities
\begin{equation}\label{eq2.1}
r^{\star}\left(f_{\alpha_d}\right)<\sqrt[d+1]{\left(d+1\right)^dA(\alpha_d)}
\end{equation}
and
 \begin{equation}\label{eq2.2}
 \frac{\left(d+1\right)^{d+1}A(\alpha_d)}{d+2}<\left(r^{\star}\left(f_{\alpha_d}\right)\right)^{d+1}<\frac{\left(d+2\right)\left(d+1\right)^{d+1}A(\alpha_d)B(\alpha_d)}{\left(d+2\right)^{2d}B(\alpha_d)-\left(2d+3\right)A(\alpha_d)}.
 \end{equation}
\end{theorem}
It is worth to mention that for $d=1$ and $\alpha_1=\nu$ the inequalities \eqref{eq2.1} and \eqref{eq2.2} reduce to the inequalities \cite[Theorem 1]{aktas1}
$$
r^{\star}\left(\varphi_\nu\right)<\sqrt{2(\nu+1)}
$$
and
$$
2\sqrt{\frac{\nu+1}{3}}< r^{\star}\left(\varphi_\nu\right)<2\sqrt{\frac{3(\nu+1)(\nu+2)}{4\nu+13}},
$$
where $\varphi_\nu$ denotes the normalized Bessel function $z\mapsto\varphi_\nu(z)=2^\nu\Gamma\left(\nu+1\right)z^{1-\nu}J_\nu(z)$.

\subsection{Radii of convexity of normalized hyper-Bessel functions}
Our second main result is on the radius of convexity $r^{c}\left(f_{\alpha_d}\right)$ of the normalized hyper-Bessel function $z\mapsto{f_{\alpha_d}(z)}$.

\begin{theorem}\label{th2}
	Let $\alpha_{i}>-1$ for $i \in\{1,2,\dots,d\}.$	Then the radius of convexity $r^{c}\left(f_{\alpha_d}\right)$ of the normalized hyper-Bessel function
$z\mapsto{f_{\alpha_d}(z)}=z\mathcal{J}_{\alpha_d}(z)$ is the smallest positive root of the equation
$$\mathcal{J}_{\alpha_d}(z)+3z\mathcal{J}_{\alpha_d}^{\prime}(z)+z^2\mathcal{J}_{\alpha_d}^{\prime\prime}(z)=0$$
and satisfies the following inequality
\begin{equation}\label{eq2.11}
	 \frac{\left(d+1\right)^{d+1}A(\alpha_d)}{\left(d+2\right)^2}<\left(r^{c}\left(f_{\alpha_d}\right)\right)^{d+1}<\frac{\left(d+1\right)^{d+1}\left(d+2\right)^2A(\alpha_d)B(\alpha_d)}{\left(d+2\right)^{4}B(\alpha_d)-\left(2d+3\right)^2A(\alpha_d)}.
	\end{equation}
\end{theorem}

Here, it is important to note that if we take $d=1$ and $\alpha_1=\nu$ in the inequality \eqref{eq2.11} we get the inequality \cite[Theorem 6]{aktas2}
$$
\frac{2\sqrt{\nu+1}}{3}< r^{c}\left(\varphi_\nu\right)<6\sqrt{\frac{(\nu+1)(\nu+2)}{56\nu+137}},
$$
where $\varphi_\nu$ denotes the normalized Bessel function $z\mapsto\varphi_\nu(z)=2^\nu\Gamma\left(\nu+1\right)z^{1-\nu}J_\nu(z)$.

\subsection{Radii of uniform convexity of normalized hyper-Bessel functions}
In this subsection we would like to present result on the radius of uniform convexity $r^{uc}\left(f_{\alpha_d}\right)$ of the normalized hyper-Bessel function $z\mapsto{f_{\alpha_d}(z)}$.

\begin{theorem}\label{ruc}
Let $\alpha_{i}>-1$ for $i \in\{1,2,\dots,d\}$. Then the radius of uniform convexity $r^{uc}\left(f_{\alpha_d}\right)$ of the normalized hyper-Bessel function
$z\mapsto{f_{\alpha_d}(z)}=z\mathcal{J}_{\alpha_d}(z)$ is the smallest positive root of the equation $$
2z^2\mathcal{J}_{\alpha_d}^{\prime\prime}(z)+5z\mathcal{J}_{\alpha_d}^{\prime}(z)+\mathcal{J}_{\alpha_d}(z)=0.$$
\end{theorem}

It is worth also to mention that if we take $d=1$ and $\alpha_1=\nu$, we reobtain a recent result on the radius of uniform convexity of Bessel functions \cite[Theorem 3.2]{Deniz} which state that, if $\nu>-1$, then the radius of uniform convexity of the function $z\mapsto2^\nu\Gamma\left(\nu+1\right)z^{1-\nu}J_\nu(z)$ is the smallest positive
root of the equation
$$
1+2r\frac{(2\nu-1)J_{\nu+1}(r)-rJ_{\nu}(r)}{J_{\nu}(r)-rJ_{\nu+1}(r)}=0.
$$

\subsection{Bounds for the first positive zeros of normalized hyper-Bessel functions}
Our next result is about the first positive zero of normalized hyper-Bessel function.
\begin{theorem}\label{th3}
	Let $\alpha_{i}>-1$ for $i\in \{1,2,\dots,d\}.$ Then, the first positive zero $j_{\alpha_d,1}$ of the hyper-Bessel function $z\mapsto \mathcal{J}_{\alpha_d}(z)$ satisfies following inequalities
\begin{equation*}\label{eq2.16}
	 \left(d+1\right)^{d+1}A(\alpha_d)<j_{\alpha_d,1}^{d+1}<\frac{\left(d+1\right)^{d+1}A(\alpha_d)B(\alpha_d)}{B(\alpha_d)-A(\alpha_d)}
\end{equation*}
and
\begin{equation*}\label{eq2.17}
	 \left(d+1\right)^{d+1}A(\alpha_d)\sqrt{\frac{B(\alpha_d)}{B(\alpha_d)-A(\alpha_d)}}<j_{\alpha_d,1}^{d+1}<
\frac{2\left(d+1\right)^{d+1}\left(A(\alpha_d)B(\alpha_d)C(\alpha_d)-\left(A(\alpha_d)\right)^2C(\alpha_d)\right)}{\left(A(\alpha_d)\right)^2-3A(\alpha_d)C(\alpha_d)+2B(\alpha_d)C(\alpha_d)},
\end{equation*}
where $C(\alpha_d)=\prod_{i=1}^{d}\left(\alpha_{i}+3\right).$
\end{theorem}
Here it is worth mentioning that by setting $d=1$ and $\alpha_1=\nu$ in Theorem \eqref{th3}, we get the inequalities \cite[5.3 and 5.4]{ismail}
$$
4(\nu+1)<j_{\nu,1}^2<4(\nu+1)(\nu+2)
$$
and
$$
4(\nu+1)\sqrt{\nu+2}<j_{\nu,1}^2<2(\nu+1)(\nu+3)
$$
for the first positive zero $j_{\nu,1}$ of the normalized classical Bessel function $J_\nu(z)$, when $\nu>-1.$

\subsection{Interlacing properties of  the zeros of normalized hyper-Bessel functions}
Our next result is related to interlacing properties of  the zeros of normalized hyper-Bessel functions.
\begin{theorem}\label{interlace}
Let $\alpha_{i}>-1$ for $i\in \{1,2,\dots,d\}.$ Then the zeros of the function $z\mapsto \mathcal{J}'_{\alpha_d}(z)$ are interlaced with those of the function $z\mapsto \mathcal{J}_{\alpha_d}(z)$.
\end{theorem}

\subsection{Redheffer-type inequalities and bounds for normalized hyper-Bessel functions}
In this subsection we find Redheffer-type inequalities, monotonicity, and bounds for normalized hyper-Bessel functions.
\begin{theorem}\label{AM}
Let $\alpha_{i}>-1$ for $i\in \{1,2,\dots,d\}.$ Then the function $h_{\alpha_d},q_{\alpha_d}:[0,j_{\alpha_d,1}^{d+1})\rightarrow [0,\infty)$ defined by
$$
h_{\alpha_d}(x)=\left(\log\left(\frac{x^{\frac{S(\alpha_d)}{d+1}}e^{-\frac{x}{(d+1)^{d+1}A(\alpha_d)}}}{J_{\alpha_d}(x^{\frac{1}{d+1}})}\right)\right)'
$$
and
$$
q_{\alpha_d}(x)=\left(\frac{x^{\frac{S(\alpha_d)}{d+1}}e^{-\frac{x}{(d+1)^{d+1}A(\alpha_d)}}}{J_{\alpha_d}(x^{\frac{1}{d+1}})}\right)
$$
are  absolutely monotonic. Here $S(\alpha_d)=\sum_{i=1}^d\alpha_i$.
\end{theorem}

Let us recall a recent result on the complete monotonicity of certain special functions \cite[Lemma 2]{zhang} which states that if an entire function $f$ with genus $0$ has only negative zeros, and $d(f)\geq m>n$, $m$, $n+1\in \mathbb{N}$ then the function ${f^{(m)}(x)}/{f^{(n)}(x)}$ is completely monotonic in $x$. Here $d(f)$ is the degree of polynomial if $f$ is a polynomial, and $d(f)=\infty$ otherwise. In view of this result and the fact
that the function $\mathcal{J}_{\mathbf{\alpha_d}}(-x)$ is of genus $0$ (see proof of the Theorem \ref{interlace}) and has only negative real zeros  for $\alpha_i>-1,$ $i\in\{1,\dots,d\}$ conform \cite[Theorem 2]{bsingh}, we obtain that ${\mathcal{J}_{\mathbf{\alpha_d}}^{(m)}(-x)}/{\mathcal{J}_{\mathbf{\alpha_d}}^{(n)}(-x)}$ is completely monotonic in $x$ for $ m>n$,  $m$, $n+1\in \mathbb{N}$.

We note that the absolutely monotonicity of $q_{\alpha_d}$ can be used to find the upper bound for the hyper-Bessel functions. In  particular, we have:

\begin{corollary}\label{cor}
Let $\alpha_{i}>-1$ for $i\in \{1,2,\dots,d\}$ and $x\in [0,j_{\alpha_d,1}^{d+1})$. Then
$$
J_{\alpha_d}(x)\leq \left( \frac{x}{d+1}\right)^{S(\alpha_d)}\frac{e^{-\frac{x^{d+1}}{(d+1)^{d+1}A(\alpha_d)}}}{A(\alpha_d)}.
$$
\end{corollary}

Finally, we obtain some bounds for real hyper-Bessel functions.

\begin{theorem}\label{RT}
If $\alpha_{i}>-1$ for $i\in \{1,2,\dots,d\}$, and $x\in (0,j_{\alpha_d,1})$ then the following sharp exponential Redheffer-type inequalities hold:
\begin{equation}\label{Redheffer}
\left(\frac{j_{\alpha_d,1}^{d+1}-x^{d+1}}{j_{\alpha_d,1}^{d+1}}\right)^{a_{\alpha_d}}\leq \mathcal{J}_{\alpha_d}(x)\leq
\left(\frac{j_{\alpha_d,1}^{d+1}-x^{d+1}}{j_{\alpha_d,1}^{d+1}}\right)^{b_{\alpha_d}},
\end{equation}
where $a_{\alpha_d}=\frac{j_{\alpha_d,1}^{d+1}}{(d+1)^{d+1}A(\alpha_d)}$ and $b_{\alpha_d}=1$ are best possible constants.
\end{theorem}

Note that for $d=1$ and $\alpha_1=\nu$ the inequality \eqref{Redheffer} reduces to  \cite[Theorem 1]{baricz_mehrez}
$$
\left(\frac{j_{\nu,1}^{2}-x^{2}}{j_{\nu,1}^{2}}\right)^{a_{\nu}}\leq \mathcal{J}_{\nu}(x)\leq
\left(\frac{j_{\nu,1}^{2}-x^{2}}{j_{\nu,1}^{2}}\right)^{b_{\nu}}, \qquad x\in (0,j_{\nu,1})
$$
with the best possible constants $a_\nu = \frac{j_{\nu,1}^2}{4(\nu+1)}$ and $b_\nu = 1$.

\section{\bf Proofs of the main results}
\setcounter{equation}{0}

\begin{proof}[\bf Proof of Theorem \ref{th1}]

First we find the radius of starlikeness of the function $z\mapsto{f_{\alpha_d}(z)}$.
For this we shall prove that for $\alpha_i>-1, i\in \{1,2,\dots,d\},$ the inequality  $\real\left(\dfrac{z{f_{\alpha_d}^{\prime}(z)}}{{f_{\alpha_d}(z)}}\right)>0$ holds for $z\in \mathbb{D}_{r^*(f_{\alpha_d})}$ and does not hold in any larger disk.  By taking the logarithmic derivative of \eqref{eq 1.10} we have
	\begin{equation}\label{eq2.3}
	 \dfrac{z{f_{\alpha_d}^{\prime}(z)}}{{f_{\alpha_d}(z)}}=1-(d+1)\sum_{n\geq1}\frac{z^{d+1}}{{j_{\alpha_d,n}^{d+1}}-z^{d+1}}.
	\end{equation}
In view of the inequality \cite{szasz}
	\begin{equation}\label{eq2.4}
	 \real\left(\frac{z}{\alpha-z}\right)\leq\frac{\left|z\right|}{\alpha-\left|z\right|},
	\end{equation}
which holds for $z\in\mathbb{C}, \alpha\in\mathbb{R}\text{ and } \left|z\right|<\alpha$, we get
	\begin{equation*}\label{eq2.5}
	 \frac{\left|z\right|^{d+1}}{{j_{\alpha_d,n}^{d+1}}-\left|z\right|^{d+1}}\geq\real\left(\frac{z^{d+1}}{{j_{\alpha_d,n}^{d+1}}-z^{d+1}}\right)
	\end{equation*}
	for all $\alpha_i>-1, n\in\{1,2,\dots\}$ and $\left|z\right|<j_{\alpha_d,1}$.
By using the above inequality we have
	\begin{align*}	 \real\left(\dfrac{z{f_{\alpha_d}^{\prime}(z)}}{{f_{\alpha_d}(z)}}\right)
=\real\left(1-(d+1)\sum_{n\geq1}\frac{z^{d+1}}{{j_{\alpha_d,n}^{d+1}}-z^{d+1}}\right)
\geq1-(d+1)\sum_{n\geq1}\frac{\left|z\right|^{d+1}}{{j_{\alpha_d,n}^{d+1}}-\left|z\right|^{d+1}}
=\dfrac{\left|z\right|{f_{\alpha_d}^{\prime}(\left|z\right|)}}{{f_{\alpha_d}(\left|z\right|)}}
	\end{align*}
and equality holds only when $z=|z|=r$.	The minimum principle for harmonic functions and the last inequality imply that
$$
\real\left(\dfrac{z{f_{\alpha_d}^{\prime}(z)}}{{f_{\alpha_d}(z)}}\right)>0 \text{ if and only if }\left|z\right|<x_{\alpha_d,1},
$$
where $x_{\alpha_d,1}$ is the smallest positive root of the equation
\begin{equation}\label{logder}
\dfrac{r{f_{\alpha_d}^{\prime}(r)}}{{f_{\alpha_d}(r)}}=0,
\end{equation}
which is equivalent to
$$
z\mathcal{J}_{\alpha_d}^{\prime}(z)+\mathcal{J}_{\alpha_d}(z)=0.
$$
	
Now, we show that the inequality \eqref{eq2.1} is true. Since $x_{\alpha_d,1}=r^{\star}\left(f_{\alpha_d}\right)$ is the smallest positive root of the equation \eqref{logder}, the equation \eqref{eq2.3}
vanish at $r^{\star}\left(f_{\alpha_d}\right).$ Consequently, we have
\begin{equation*}\label{eq2.6}
\frac{1}{(d+1)\left(r^{\star}\left(f_{\alpha_d}\right)\right)^{d+1}}=\sum_{n\geq1}\frac{1}
{{j_{\alpha_d,n}^{d+1}}-\left(r^{\star}\left(f_{\alpha_d}\right)\right)^{d+1}}>\sum_{n\geq1}\frac{1}{{j_{\alpha_d,n}^{d+1}}},
\end{equation*}
which in view of the equation \eqref{eq.1.13} gives the inequality \eqref{eq2.1}.

The inequality \eqref{eq2.2} can be obtained by using the Euler-Rayleigh inequalities and some properties of Laguerre-P\'{o}lya class $\mathcal{LP}$ of real entire functions. Note that all zeros of the hyper-Bessel function $z\mapsto \mathcal{J}_{\alpha_d}(z)$ are real for $\alpha_{i}>-1$ and $i\in\{1,2,\dots,d\}$ (see \cite[Theorem 2]{bsingh}). Thus, the function
$$
z\mapsto{f_{\alpha_d}(z)}=z\mathcal{J}_{\alpha_d}(z)
$$ has also real zeros for $\alpha_{i}>-1$ and $i\in\{1,2,\dots,d\}$. Therefore, the function $z\mapsto{f_{\alpha_d}(z)}$ belongs to the Laguerre-P\'{o}lya class $\mathcal{LP}$ of real entire functions. Since the class $\mathcal{LP}$ is closed under differentiation, the function
\begin{equation}\label{eq2.7}
z\mapsto\Psi_{\alpha_d}(z)={f_{\alpha_d}^{\prime}(z)}=\sum_{n\geq0}\frac{\left(n(d+1)+1\right)(-1)^n}{n!\left(d+1\right)^{n(d+1)}(\alpha_1+1)_n\cdots(\alpha_d+1)_n}z^{n(d+1)}
\end{equation}
belongs also to the  Laguerre-P\'{o}lya class $\mathcal{LP}$ of real entire functions. On the other hand, the order of an entire function $a_0+a_1z+a_2z^2+\cdots$ can be determined by the formula \cite[p. 6]{lev}
$$
\rho=\limsup_{n\to\infty}\frac{n\log{n}}{-\log\left|a_n\right|}.
$$
Since
$$
\limsup_{n\to\infty}\frac{n\log{n}}{\log{n!}+n(d+1)\log(d+1)+\sum_{i=1}^{d}\log{(\alpha_i+1)_n}-\log(n(d+1)+1)}=\frac{1}{1+d},
$$
the growth order of the function $z\mapsto\Psi_{\alpha_d}(z)$ is $\rho=\frac{1}{1+d}<1.$ Therefore by applying Hadamard's theorem \cite[p. 26]{lev} the function $z\mapsto\Psi_{\alpha_d}(z)$ can be written by the canonical product as follows:
\begin{equation}\label{eq2.8}
\Psi_{\alpha_d}(z)=\prod_{n\geq1}\left(1-\frac{z^{d+1}}{\psi_{\alpha_d,n}^{d+1}}\right),
\end{equation}
where ${\psi_{\alpha_d,n}}$ denotes the $n$th positive real zero of the function $\Psi_{\alpha_d}.$ Now, taking logarithmic derivative of both sides of \eqref{eq2.8} we have
\begin{equation}\label{eq2.9}
\frac{\Psi_{\alpha_d}^{\prime}(z)}{\Psi_{\alpha_d}(z)}=-(d+1)\sum_{k\geq0}\delta_{k+1}z^{k(d+1)+d},\qquad   \left|z\right|<{\psi_{\alpha_d,1}},
\end{equation}
where
$$
\delta_k=\sum_{n\geq1}\left({\psi_{\alpha_d,n}}\right)^{-k(d+1)},
$$
is Euler-Rayleigh sum for the zeros of $\Psi_{\alpha_d}.$  Also, using the infinite sum representation \eqref{eq2.7} we get
\begin{equation}\label{eq2.10}
\frac{\Psi_{\alpha_d}^{\prime}(z)}{\Psi_{\alpha_d}(z)}={\sum_{n\geq0}U_nz^{(n+1)(d+1)-1}}\Bigg{/}{\sum_{n\geq0}V_nz^{n(d+1)}},
\end{equation} where
$$
U_n=\frac{(-1)^{n+1}(n+1)(d+1)\left((n+1)(d+1)+1\right)}{(n+1)!\left(d+1\right)^{(n+1)(d+1)}(\alpha_1+1)_{n+1}\cdots(\alpha_d+1)_{n+1}}
$$
and
$$
V_n=\frac{(-1)^{n}\left(n(d+1)+1\right)}{n!\left(d+1\right)^{n(d+1)}(\alpha_1+1)_{n}\cdots(\alpha_d+1)_{n}}.
$$
By comparing the coefficients with the same degrees of \eqref{eq2.9} and \eqref{eq2.10} we obtain the Euler-Rayleigh sums
$$
\delta_1=\frac{d+2}{(d+1)^{d+1}A(\alpha_d)}
$$
and
$$
\delta_2=\frac{(d+2)^2B(\alpha_d)-(2d+3)A(\alpha_d)}{(d+1)^{2(d+1)}\left(A(\alpha_d)\right)^2B(\alpha_d)}.
$$
By using the Euler-Rayleigh inequalities
$$
\delta_k^{-\frac{1}{k}}<{\psi_{\alpha_d,1}^{d+1}}<\frac{\delta_k}{\delta_{k+1}}
$$
for $\alpha_{i}>-1$, $i\in\{1,2,\dots,d\}$ and $k=1$ we get the following inequality
$$
\frac{\left(d+1\right)^{d+1}A(\alpha_d)}{d+2}
<\left(r^{\star}\left(f_{\alpha_d}(z)\right)\right)^{d+1}
<\frac{\left(d+2\right)\left(d+1\right)^{d+1}A(\alpha_d)B(\alpha_d)}{\left(d+2\right)^{2}B(\alpha_d)-\left(2d+3\right)A(\alpha_d)}.
$$

It is possible to obtain more tighter lower and upper bounds for the radius of starlikeness of the function $z\mapsto f_{\alpha_d}(z)$, but we restricted ourselves for first two values of $k$ since the other values would be quite complicated.
\end{proof}

\begin{proof}[\bf Proof of Theorem \ref{th2}]
Let us recall the Alexander's duality theorem \cite{Alexander} which states that, if $f:\mathbb{D}\rightarrow \mathbb{C}$ be an analytic function then $f$ is
convex if and only if $zf'$ is starlike. Therefore it is enough to find the radius of starlikeness of  the function $zf_{\alpha_d}^{\prime}$.
Now let us consider $F_{\alpha_d}(z)=zf_{\alpha_d}^{\prime}(z)$, which in view of \eqref{eq2.7} and \eqref{eq2.8} gives
$$
F_{\alpha_d}(z)=z\Psi_{\alpha_d}(z)=z\prod_{n\geq1}\left(1-\frac{z^{d+1}}{\psi_{\alpha_d,n}^{d+1}}\right).
$$
By taking the logarithmic derivative of the above equation we have
	\begin{equation*}
	 \dfrac{z{F_{\alpha_d}^{\prime}(z)}}{{F_{\alpha_d}(z)}}=1-(d+1)\sum_{n\geq1}\frac{z^{d+1}}{{\psi_{\alpha_d,n}^{d+1}}-z^{d+1}}.
	\end{equation*}
Now in view of the inequality \eqref{eq2.4} we get
	\begin{equation*}\label{eq2.5}
	 \frac{\left|z\right|^{d+1}}{{\psi_{\alpha_d,n}^{d+1}}-\left|z\right|^{d+1}}\geq\real\left(\frac{z^{d+1}}{{\psi_{\alpha_d,n}^{d+1}}-z^{d+1}}\right)
	\end{equation*}
	for all $\alpha_i>-1, n\in\{1,2,\dots\}$ and $\left|z\right|<\psi_{\alpha_d,1}$.
By using the above inequality we have
\begin{align*}
\real\left(\dfrac{z{F_{\alpha_d}^{\prime}(z)}}{{F_{\alpha_d}(z)}}\right)&=\real\left(1-(d+1)\sum_{n\geq1}\frac{z^{d+1}}
{{\psi_{\alpha_d,n}^{d+1}}-z^{d+1}}\right)
\geq 1-(d+1)\sum_{n\geq1}\frac{\left|z\right|^{d+1}}{{\psi_{\alpha_d,n}^{d+1}}-\left|z\right|^{d+1}}
=\dfrac{\left|z\right|{F_{\alpha_d}^{\prime}(\left|z\right|)}}{{F_{\alpha_d}(\left|z\right|)}}
	\end{align*}
and equality holds only when $z=|z|=r$.	
 The minimum principle for harmonic functions and the last inequality imply that
$$
\real\left(\dfrac{z{F_{\alpha_d}^{\prime}(z)}}{{F_{\alpha_d}(z)}}\right)>0 \text{ if and only if }\left|z\right|<y_{\alpha_d,1},
$$
where $y_{\alpha_d,1}$ is the smallest positive root of the equation
\begin{equation*}\label{logder2}
\dfrac{r{F_{\alpha_d}^{\prime}(r)}}{{F_{\alpha_d}(r)}}=0,
\end{equation*}
which is equivalent to
$$
\mathcal{J}_{\alpha_d}(z)+3z\mathcal{J}_{\alpha_d}^{\prime}(z)+z^2\mathcal{J}_{\alpha_d}^{\prime\prime}(z)=0.
$$
This proves the first part of Theorem \ref{th2}.
	
Now, we prove that the radius of convexity $r^{c}\left(f_{\alpha_d}\right)$ of the function $z\mapsto{f_{\alpha_d}(z)}$ satisfies the inequality \eqref{eq2.11}. By using the infinite sum representation \eqref{eq2.7} we have
	\begin{equation}\label{eq2.12}
	 \Phi_{\alpha_d}(z)=\left(zf_{\alpha_d}^{\prime}(z)\right)^{\prime}=\sum_{n\geq0}\frac{(-1)^n\left(n(d+1)+1\right)^2}{n!\left(d+1\right)^{n(d+1)}(\alpha_1+1)_n\cdots(\alpha_d+1)_n}z^{n(d+1)}.
	\end{equation}
	We know that the function $z\mapsto{f_{\alpha_d}(z)}$ belongs to the Laguerre-P\'{o}lya class $\mathcal{LP}$ of real entire functions. Considering the properties of the Laguerre-P\'{o}lya class $\mathcal{LP}$ of real entire functions we can say that the function $z\mapsto{\Phi_{\alpha_d}(z)}$ also belongs  to the Laguerre-P\'{o}lya class $\mathcal{LP}$ of real entire functions. Since
$$
\rho=\limsup_{n\to\infty}\frac{n\log{n}}{\log{n!}+n(d+1)\log(d+1)+\sum_{i=1}^{d}\log{(\alpha_i+1)_n}-2\log(n(d+1)+1)}=\frac{1}{d+1},
$$
according to Hadamard's theorem \cite[p. 26]{lev}, the function $z\mapsto{\Phi_{\alpha_d}(z)}$ can be written by the infinite product
	\begin{equation}\label{eq2.13}
	 \Phi_{\alpha_d}(z)=\prod_{n\geq1}\left(1-\frac{z^{d+1}}{\tau_{\alpha_d,n}^{d+1}}\right),
	\end{equation}
	where $\tau_{\alpha_d,n}$ denotes the $n$th positive zero of the function $\Phi_{\alpha_d}$. Now, the Euler-Rayleigh sums can be expressed in terms of $d$ and $\alpha_i$ for $i\in\{1,2,\dots,d\}$. By taking logarithmic derivative of \eqref{eq2.13} we get
	\begin{equation}\label{eq2.14}
	 \frac{\Phi_{\alpha_d}^{\prime}(z)}{\Phi_{\alpha_d}(z)}=-(d+1)\sum_{k\geq0}\epsilon_{k+1}z^{k(d+1)+d}, \qquad  \left|z\right|<{\tau_{\alpha_d,1}},
	\end{equation}
	where
$$
\epsilon_k=\sum_{n\geq1}\left({\tau_{\alpha_d,n}}\right)^{-k(d+1)}
$$
is Euler-Rayleigh sum for the zeros of $\Phi_{\alpha_d}.$  Also, using the infinite sum representation \eqref{eq2.12} we get
	\begin{equation}\label{eq2.15}
	 \frac{\Phi_{\alpha_d}^{\prime}(z)}{\Phi_{\alpha_d}(z)}={\sum_{n\geq0}S_nz^{n(d+1)+d}}\Bigg{/}{\sum_{n\geq0}T_nz^{n(d+1)}},
	\end{equation}
	where
$$
S_n=\frac{(-1)^{n+1}\left((n+1)(d+1)+1\right)^2(n+1)(d+1)}{(n+1)!\left(d+1\right)^{(n+1)(d+1)}(\alpha_1+1)_{n+1}\cdots(\alpha_d+1)_{n+1}}
$$
and
$$
T_n=\frac{(-1)^{n}\left(n(d+1)+1\right)^2}{n!\left(d+1\right)^{n(d+1)}(\alpha_1+1)_{n}\cdots(\alpha_d+1)_{n}}.
$$
By equating the right-hand sides of the equations \eqref{eq2.14} and \eqref{eq2.15} we have the following Euler-Rayleigh sums
$$
\epsilon_1=\frac{(d+2)^2}{(d+1)^{d+1}A(\alpha_d)}
$$
and
$$
\epsilon_2=\frac{(d+2)^4B(\alpha_d)-(2d+3)^2A(\alpha_d)}{(d+1)^{2(d+1)}\left(A(\alpha_d)\right)^2B(\alpha_d)}.
$$
If we consider these Euler-Rayleigh sums in the Euler-Rayleigh inequalities
$$
\epsilon_k^{-\frac{1}{k}}<{\tau_{\alpha_d,1}^{d+1}}<\frac{\epsilon_k}{\epsilon_{k+1}}
$$ we obtain the inequality \eqref{eq2.11}.
\end{proof}

\begin{proof}[\bf Proof of Theorem \ref{ruc}]
In view of \eqref{eq2.7} and \eqref{eq2.8}, we have
\begin{equation}\label{1.18}
1+\frac{z{f_{\alpha_d}^{\prime\prime}(z)}}{{f_{\alpha_d}^{\prime}(z)}}=1-(d+1)\sum_{n\geq1}\frac{z^{d+1}}{\psi_{\alpha_d,n}^{d+1}-z^{d+1}}.
\end{equation}

By using the equality \eqref{1.18} and the inequality \eqref{eq2.4}  we have that
\begin{align*}
\real\left(1+\frac{z{f_{\alpha_d}^{\prime\prime}(z)}}{{f_{\alpha_d}^{\prime}(z)}}\right)
&=1-\real\left((d+1)\sum_{n\geq1}\frac{z^{d+1}}{\psi_{\alpha_d,n}^{d+1}-z^{d+1}}\right)\\&\geq1-(d+1)\left(\sum_{n\geq1}\frac{r^{d+1}}{\psi_{\alpha_d,n}^{d+1}-r^{d+1}}\right)
=1+r\frac{{f_{\alpha_d}^{\prime\prime}(r)}}{{f_{\alpha_d}^{\prime}(r)}},
\end{align*}
where $\left|z\right|\leq{r}<\psi_{\alpha_d,1}.$ Therefore, we get that
\begin{equation}\label{2.1}
\real\left(1+\frac{z{f_{\alpha_d}^{\prime\prime}(z)}}{{f_{\alpha_d}^{\prime}(z)}}\right)
\geq1+r\frac{{f_{\alpha_d}^{\prime\prime}(r)}}{{f_{\alpha_d}^{\prime}(r)}}, \qquad \left|z\right|\leq{r}<\psi_{\alpha_d,1}.
\end{equation}
On the other hand, the equation \eqref{1.18} implies that
\begin{align*}
\left|\frac{z{f_{\alpha_d}^{\prime\prime}(z)}}{{f_{\alpha_d}^{\prime}(z)}}\right|
=\left|-(d+1)\sum_{n\geq1}\frac{z^{d+1}}{\psi_{\alpha_d,n}^{d+1}-z^{d+1}}\right|
\leq\sum_{n\geq1}\left|\frac{(d+1)z^{d+1}}{\psi_{\alpha_d,n}^{d+1}-z^{d+1}}\right|
\leq\sum_{n\geq1}\frac{(d+1)r^{d+1}}{\psi_{\alpha_d,n}^{d+1}-r^{d+1}}
=-r\frac{{f_{\alpha_d}^{\prime\prime}(r)}}{{f_{\alpha_d}^{\prime}(r)}}.
\end{align*}
Thus, we have
\begin{equation}\label{2.2}
\left|\frac{z{f_{\alpha_d}^{\prime\prime}(z)}}{{f_{\alpha_d}^{\prime}(z)}}\right|
\leq-r\frac{{f_{\alpha_d}^{\prime\prime}(r)}}{{f_{\alpha_d}^{\prime}(r)}}, \qquad \left|z\right|\leq{r}<\psi_{\alpha_d,1}.
\end{equation}
The inequalities \eqref{2.1} and \eqref{2.2} imply that
\begin{equation}\label{2.3}
\real\left(1+\frac{z{f_{\alpha_d}^{\prime\prime}(z)}}{{f_{\alpha_d}^{\prime}(z)}}\right)
-\left|\frac{z{f_{\alpha_d}^{\prime\prime}(z)}}{{f_{\alpha_d}^{\prime}(z)}}\right|
\geq1+2r\frac{{f_{\alpha_d}^{\prime\prime}(r)}}{{f_{\alpha_d}^{\prime}(r)}}, \qquad \left|z\right|\leq{r}<\psi_{\alpha_d,1}
\end{equation}
and equality  holds in \eqref{2.3} if and only if $z=r$. Thus, it follows that
$$
\inf_{\left|z\right|<r}\left(\real\left(1+\frac{z{f_{\alpha_d}^{\prime\prime}(z)}}{{f_{\alpha_d}^{\prime}(z)}}\right)
-\left|\frac{z{f_{\alpha_d}^{\prime\prime}(z)}}{{f_{\alpha_d}^{\prime}(z)}}\right|\right)
=1+2r\frac{{f_{\alpha_d}^{\prime\prime}(r)}}{{f_{\alpha_d}^{\prime}(r)}}, \qquad r\in(0,\psi_{\alpha_d,1}).
$$
The mapping $\Theta_{\alpha_d}(r):(0,\psi_{\alpha_d,1})\mapsto\mathbb{R},$ defined by
$$
\Theta_{\alpha_d}(r)=1+2r\frac{{f_{\alpha_d}^{\prime\prime}(r)}}{{f_{\alpha_d}^{\prime}(r)}}
=1-2(d+1)\sum_{n\geq1}\frac{r^{d+1}}{\psi_{\alpha_d,n}^{d+1}-r^{d+1}},
$$
is strictly decreasing, $\lim_{r\searrow0}\Theta_{\alpha_d}(r)=1$ and $\lim_{r\nearrow\psi_{\alpha_d,1}}\Theta_{\alpha_d}(r)=-\infty.$ Consequently, the equation
\begin{equation*}\label{2.4}
\Theta_{\alpha_d}(r)=1+2r\frac{{f_{\alpha_d}^{\prime\prime}(r)}}{{f_{\alpha_d}^{\prime}(r)}}=0
\end{equation*}
has a unique root $r_0$ in the interval $(0,\psi_{\alpha_d,1})$ and $r_0=r^{uc}(f_{\alpha_d})$.  Therefore, in view of the equation ${f_{\alpha_d}(z)}=z\mathcal{J}_{\alpha_d}(z)$ we obtain that the radius of uniform convexity $r^{uc}(f_{\alpha_d})$ of the normalized hyper-Bessel function $z\mapsto{f_{\alpha_d}}(z)$ is the smallest positive root of the equation
$$
2z^2\mathcal{J}_{\alpha_d}^{\prime\prime}(z)+5z\mathcal{J}_{\alpha_d}^{\prime}(z)+\mathcal{J}_{\alpha_d}(z)=0.
$$
\end{proof}

\begin{proof}[\bf Proof of Theorem \ref{th3}]

	By using the infinite sum representation \eqref{eq1.6} we have
	\begin{equation}\label{eq2.18} \mathcal{J}_{\alpha_d}^{\prime}(z)=\sum_{n\geq0}\frac{(-1)^{n+1}(n+1)}{(n+1)!(\alpha_1+1)_{n+1}\cdots(\alpha_d+1)_{n+1}}\left(\frac{z}{d+1}\right)^{nd+n+d}.
	\end{equation}
	From the equations \eqref{eq1.6} and \eqref{eq2.18} we get
	\begin{equation}\label{eq2.19}
	 \frac{\mathcal{J}_{\alpha_d}^{\prime}(z)}{\mathcal{J}_{\alpha_d}(z)}=\sum_{n\geq0}\kappa_nz^{nd+n+d}\Bigg{/}\sum_{n\geq0}\xi_nz^{nd+n},
	\end{equation}
	where
$$
\kappa_n=\frac{(-1)^{n+1}(n+1)}{(n+1)!(d+1)^{nd+n+d}(\alpha_1+1)_{n+1}\cdots(\alpha_d+1)_{n+1}}
$$
and
$$
\xi_n=\frac{(-1)^{n}}{n!(d+1)^{nd+n}(\alpha_1+1)_{n}\cdots(\alpha_d+1)_{n}}.
$$
Also, logarithmic differentiation of the equation \eqref{eq 1.9} implies that
	\begin{equation}\label{eq2.20}
	 \frac{\mathcal{J}_{\alpha_d}^{\prime}(z)}{\mathcal{J}_{\alpha_d}(z)}=-(d+1)\sum_{k\geq0}\Delta_{k+1}z^{kd+k+d}, \qquad  \left|z\right|<{j_{\alpha_d,1}},
	\end{equation}
where
$$
\Delta_{k}=\sum_{n\geq1}j_{\alpha_d,n}^{-k(d+1)}
$$
is the Euler-Rayleigh sum for the zeros of $\mathcal{J}_{\alpha_d}.$ Now, it is possible to express the Euler-Rayleigh sums $\Delta_{k}$ in terms of $d$ and $\alpha_i$ for $i\in\{1,2,\dots,d\}.$ Comparing the coefficients of the equations \eqref{eq2.19} and \eqref{eq2.20} we obtain the following Euler-Rayleigh sums, that is,
$$
\Delta_1=\frac{1}{(d+1)^{d+1}A(\alpha_d)}, \Delta_2=\frac{B(\alpha_d)-A(\alpha_d)}{(d+1)^{2(d+1)}\left(A(\alpha_d)\right)^2B(\alpha_d) }
$$
and
$$
\Delta_3=\frac{\left(A(\alpha_d)\right)^2-3A(\alpha_d)C(\alpha_d)+2B(\alpha_d)C(\alpha_d)}
{2(d+1)^{3(d+1)}\left(A(\alpha_d)\right)^3B(\alpha_d)C(\alpha_d)}.
$$
Now, using the Euler-Rayleigh inequalities
$$
\Delta_{k}^{-\frac{1}{k}}<j_{\alpha_d,1}^{d+1}<\frac{\Delta_{k}}{\Delta_{k+1}}
$$
for $k\in\{1,2\}$, we get the desired inequalities.
\end{proof}

\begin{proof}[\bf Proof of Theorem \ref{interlace}]

The function $z\mapsto \mathcal{J}_{\alpha_d}(z)$ is an entire functions of order $\frac{1}{d+1}$. Therefore
the genus of the entire function $z\mapsto \mathcal{J}_{\alpha_d}(z)$ is $0$, as the genus of entire function of order $\rho$ is
$[\rho]$ when $\rho$ is not an integer \cite[p. 34]{boas}. We also note that the zeros of $z\mapsto \mathcal{J}_{\alpha_d}(z)$ are all real
when $\alpha_{i}>-1$ for $i\in \{1,2,\dots,d\}.$ Now recall Laguerre's theorem on separation of zeros \cite[p. 23]{boas} which states that, if $z\mapsto f (z)$ is a
non-constant entire function, which is real for real $z$ and has only real zeros, and is of genus $0$ or $1$, then the zeros
of $f'$ are also real and separated by the zeros of $f$ . Therefore in view of Laguerre's theorem the conclusions follow.
\end{proof}

\begin{proof}[\bf Proof of Theorem \ref{AM}]
We note that
$$
h_{\alpha_d}(x)=\left(\log\left(\frac{x^{\frac{S(\alpha_d)}{d+1}}e^{-\frac{x}{(d+1)^{d+1}A(\alpha_d)}}}{J_{\alpha_d}(x^{\frac{1}{d+1}})}\right)\right)'=\sum_{n\geq 1}\frac{1}{j_{\alpha_d,n}^{d+1}-x}-\frac{1}{(d+1)^{d+1}A(\alpha_d)}.
$$
Therefore by differentiating $m$ times we obtain
$$
h_{\alpha_d}^{(m)}(x)=\sum_{n\geq 1}\frac{m!}{(j_{\alpha_d,n}^{d+1}-x)^m}\geq 0
$$
for all $m\in \mathbb{N}$. Hence $x\mapsto h_{\alpha_d}(x)$ is increasing on  $[0,j_{\alpha_d,1}^{d+1})$ for all $\alpha_i>-1$ and in view of the equation \eqref{eq.1.13}, $h_{\alpha_d}(0)=0$. Therefore, $h_{\alpha_d}(x)\geq h_{\alpha_d}(0)=0$. This proves the absolutely monotonicity of $x\mapsto h_{\alpha_d}(x)$ on  $[0,j_{\alpha_d,1}^{d+1})$ for all $\alpha_i>-1$, $i\in \{1,2,\dots,d\}.$

Now, by using the fact that the exponential of a function having an absolutely monotonic derivative is absolutely monotonic, we conclude that $x\mapsto q_{\alpha_d}(x)$ is absolutely monotonic on  $[0,j_{\alpha_d,1}^{d+1})$ for all $\alpha_i>-1$, $i\in \{1,2,\dots,d\}.$
\end{proof}

\begin{proof}[\bf Proof of Corollary \ref{cor}]
Since $x\mapsto q_{\alpha_d}(x)$ is
absolutely monotonic on $[0,j_{\alpha_d,1}^{d+1})$, it is increasing on $[0,j_{\alpha_d,1}^{d+1})$. Therefore we get
$$
q_{\alpha_d}(x)\geq q_{\alpha_d}(0)=(d+1)^{S(\alpha_d)}A(\alpha_d),
$$
which implies that
$$
J_{\alpha_d}(x^{\frac{1}{d+1}})\leq  \frac{x^{\frac{S(\alpha_d)}{d+1}}e^{-\frac{x}{(d+1)^{d+1}A(\alpha_d)}}}{(d+1)^{S(\alpha_d)}A(\alpha_d)}.
$$
Hence by changing $x$ to $x^{d+1}$ we get the required inequality.
\end{proof}

\begin{proof}[\bf Proof of Theorem \ref{RT}]
We define a function $\Sigma_{\alpha_d}:(0,j_{\alpha_d,1})\rightarrow \mathbb{R}$ by
$$
\Sigma_{\alpha_d}(x)=\frac{\log\mathcal{J}_{\alpha_d}(x)}{\log\left(\frac{j_{\alpha_d,1}^{d+1}-x^{d+1}}{j_{\alpha_d,1}^{d+1}}\right)}.
$$
By virtue of \eqref{eq 1.9} note that
 $\mathcal{J}_{\alpha_d}(x)>0$ for $x\in (0,j_{\alpha_d,1})$,  $\alpha_{i}>-1$, $i\in \{1,2,\dots,d\}$. Hence  $\Sigma_{\alpha_d}$ is well defined. In view of \eqref{eq2.20} we have
 \begin{eqnarray*}
\frac{\frac{d}{dx}\log\mathcal{J}_{\alpha_d}(x)}{\frac{d}{dx}\log\left(\frac{j_{\alpha_d,1}^{d+1}-x^{d+1}}{j_{\alpha_d,1}^{d+1}}\right)}
=\frac{\sum_{m\geq 1}\Delta_{m}x^{m(d+1)}}{\sum_{m\geq 1}j_{\alpha_d,1}^{-m(d+1)}x^{m(d+1)}},
 \end{eqnarray*}
where
$$
\Delta_{m}=\sum_{n\geq1}j_{\alpha_d,n}^{-m(d+1)}
$$
is Euler-Rayleigh sum for the zeros of $\mathcal{J}_{\alpha_d}.$
Consider the sequence $\{s_m\}_{m\geq 1}$ defined by
$$
s_m=j_{\alpha_d,1}^{m(d+1)}\Delta_{m}.
$$
Since
\begin{eqnarray*}
s_{m+1}-s_m&=&j_{\alpha_d,1}^{(m+1)(d+1)}\Delta_{m+1}-j_{\alpha_d,1}^{m(d+1)}\Delta_{m}
=j_{\alpha_d,1}^{m(d+1)}\left(\sum_{n\geq 1}\frac{1}{j_{\alpha_d,n}^{m(d+1)}}\left(\frac{j_{\alpha_d,1}^{d+1}}{j_{\alpha_d,n}^{d+1}}-1\right)\right)<0
\end{eqnarray*}
for all $m\in \mathbb{N}$. Hence $\{s_m\}_{m\geq 1}$ is a decreasing sequence. Now making use of monotonicity result for power series \cite[Lemma 2.1]{pv}, the
function
$$
x \mapsto \frac{\frac{d}{dx}\log\mathcal{J}_{\alpha_d}(x)}{\frac{d}{dx}\log\left(\frac{j_{\alpha_d,1}^{d+1}-x^{d+1}}{j_{\alpha_d,1}^{d+1}}\right)}
$$
is decreasing on $(0,j_{\alpha_d,1})$ and consequently the monotone form of l'Hospital's rule \cite[Lemma 2.2]{anderson} implies that  $\Sigma_{\alpha_d}$ is decreasing on $(0,j_{\alpha_d,1})$ and hence
$$
\lim_{x\nearrow j_{\alpha_d,1}}\Sigma_{\alpha_d}(x)< \Sigma_{\alpha_d}(x) < \lim_{x\searrow 0}\Sigma_{\alpha_d}(x),
$$
which in view of the limits
$$
\lim_{x\searrow 0}\Sigma_{\alpha_d}(x)= j_{\alpha_d,1}^{d+1}\Delta_{1}=\frac{j_{\alpha_d,1}^{d+1}}{(d+1)^{d+1}A(\alpha_d)}=a_{\alpha_d}\ \mbox{and}\
\lim_{x\nearrow j_{\alpha_d,1}}\Sigma_{\alpha_d}(x)=1=b_{\alpha_d}
$$
gives the inequality \eqref{Redheffer}.
\end{proof}


\begin{thebibliography}{width}

\bibitem{aktas3}
	\textsc{\.{I}. Akta\c{s}, \'A. Baricz}, Bounds for radii of starlikeness of some $q$-Bessel functions, {\em Results Math.} 72(1) (2017) 947--963.
	
\bibitem{aktas2}
	\textsc{\.{I}. Akta\c{s}, \'A. Baricz, H. Orhan}, Bounds for the radii of starlikeness and convexity of some special functions, {\em Turkish J. Math.} 42(1) (2018) 211--226.
	
	
\bibitem{aktas1}
	\textsc{\.{I}. Akta\c{s}, \'A. Baricz, N. Ya\u{g}mur}, Bounds for the radii of univalence of some special functions, {\em Math. Inequal. Appl.} 20(3) (2017) 825--843.
	
\bibitem{aktas4}
	\textsc{\.{I}. Akta\c{s}, H. Orhan}, Bounds for the radii of convexity of some $q$-Bessel functions, arXiv:1702.04549.
	
\bibitem{Alexander}
	\textsc{J. W. Alexander}, Functions which map the interior of the unit circle upon simple regions, {\em  Ann. of Math.} 17(1) (1915) 12--22.

\bibitem{anderson}
\textsc{G.D. Anderson, M.K. Vamanamurthy, M. Vuorinen }, Inequalities for quasiconformal mappings in space, {\em Pacific J. Math.} 160(1) (1993) 1--18.
	
\bibitem{mathematica} \textsc{\'A. Baricz}, Geometric properties of generalized Bessel functions of complex order, {\em Mathematica} 48(71) (2006) 13--18.
	
\bibitem{publ} \textsc{\'A. Baricz}, Geometric properties of generalized Bessel functions, {\em Publ. Math. Debrecen} 73 (2008) 155--178.
	
\bibitem{lecture} \textsc{\'A. Baricz}, {\em Generalized Bessel Functions of the First Kind}, Lecture Notes in Mathematics, vol. 1994, Springer-Verlag, Berlin, 2010.
	
\bibitem{bdoy} \textsc{\'A. Baricz, D.K. Dimitrov, H. Orhan, N. Ya\u{g}mur}, Radii of starlikeness of some special functions, {\em Proc. Amer. Math. Soc.} 144(8) (2016) 3355--3367.
	
\bibitem{bsk} \textsc{\'A. Baricz, P.A. Kup\'an, R. Sz\'asz}, The radius of starlikeness of normalized Bessel functions of the first kind, {\em Proc. Amer. Math. Soc.} 142(6) (2014) 2019--2025.

\bibitem{baricz_mehrez}
	\textsc{\'A. Baricz, K. Mehrez}, Redheffer type bounds for Bessel and modified Bessel functions of the first kind, {\em Aequationes Math.} (accepted.) arXiv:1701.08446.

\bibitem{samy} \textsc{\'A. Baricz, S. Ponnusamy}, Starlikeness and convexity of generalized Bessel functions, {\em Integr. Transforms Spec. Funct.} 21 (2010) 641--653.	
	
\bibitem{bsingh} \textsc{\'A. Baricz, S. Singh}, Zeros of some special entire functions, {\em Proc. Amer. Math. Soc.} (in press).
	
\bibitem{btk} \textsc{\'A. Baricz, E. Toklu, E. Kadio\u{g}lu}, Radii of starlikeness and convexity of Wright functions, {\em  Math. Commun.} 23(1) (2018) 97--117.
	
\bibitem{bos} \textsc{\'A. Baricz, H. Orhan, R. Sz\'asz}, The radius of $\alpha$-convexity of normalized Bessel functions of the first kind, {\em Comput. Methods Funct. Theory} 16(1) (2016) 93--103.
	
	
\bibitem{basz} \textsc{\'A. Baricz, R. Sz\'asz}, The radius of convexity of normalized Bessel functions of the first kind, {\em Anal. Appl.} 12(5) (2014) 485--509.
	
\bibitem{boas}
\textsc{R.P. Boas}, {\em Entire functions}, Academic Press Inc., New York, 1954.	
	
\bibitem{Chaggara} \textsc{H. Chaggara, N. B. Romdhane}, On the zeros of the hyper-Bessel function, {\em Integral Transforms Spec. Funct.} 26(2) (2015) 96--101.

\bibitem{Deniz}
\textsc{E. Deniz, R. Sz\'asz}, The radius of uniform convexity of Bessel functions, {\em  J. Math. Anal. Appl.} 453(1) (2017) 572--588.
	
	
\bibitem{Duren}
\textsc{P.L. Duren}, {\em Univalent Functions}, Grundlehren Math. Wiss. 259, Springer, New York, 1983.

\bibitem{Goodman}
\textsc{A.W. Goodman}, On uniformly convex functions, {\em Ann. Polon. Math.} 56 (1991) 87--92.
		
\bibitem{ismail} \textsc{M.E.H. Ismail, M.E. Muldoon}, Bounds for the small real and purely imaginary zeros of Bessel and related functions,
	{\em Methods Appl. Anal.} 2(1) (1995) 1--21.
	
\bibitem{lev} \textsc{B.Ya. Levin}, {\em Lectures on Entire Functions}, Transl. of Math. Monographs, vol. 150, Amer. Math. Soc., 1996.

\bibitem{pv}
\textsc{S. Ponnusamy, M. Vuorinen}, Asymptotic expansions and inequalities for hypergeometric functions,
{\em Mathematika} 44 (1997), 278--301.

\bibitem{Ronning}
\textsc{F. R\o{}nning}, Uniformly convex functions and a corresponding class of starlike functions,	{\em Proc. Amer. Math. Soc.} 118(1) (1993) 189--196.
		
\bibitem{szasz} \textsc{R. Sz\'asz}, On starlikeness of Bessel functions of the first kind, In: Proceedings of the 8th Joint Conference on Mathematics	and Computer Science, Kom\'arno, Slovakia, 2010, 9pp.

\bibitem{zhang} \textsc{R. Zhang}, On complete monotonicity of certain special functions, {\em Proc. Amer. Math. Soc.} (in press).

\end{thebibliography}
\end{document}